\documentstyle[12pt]{article}
\newtheorem{thm}{Theorem}[section]
\newtheorem{dfn}{Definition}[section]
\newtheorem{prop}{Proposition}[section]
\newtheorem{lem}{Lemma}[section]
\newtheorem{rem}{Remark}[section]
\newtheorem{cor}{Corollary}[section]

\newtheorem{conj}{Conjecture}[section]
\newtheorem{ex}{Example}[section]

\begin{document} 
\title{Exterior differential algebras and flat connections on 
Weyl groups} 
\date{}
\author{Anatol N. Kirillov and Toshiaki Maeno}
\maketitle
\footnote{Both of the authors were 
supported by Grant-in-Aid for Scientific Research.}
\begin{abstract}
We study some aspects of noncommutative differential geometry 
on a finite Weyl group in the sense of S. Woronowicz, 
K. Bresser {\it et al.}, and S. Majid. For any finite Weyl 
group $W$ we consider the subalgebra generated 
by flat connections in the left-invariant exterior differential 
algebra of $W.$ 

For root systems of type $A$ and $D$ we describe 
a set of relations between the flat connections, which 
conjecturally is a complete set. 
\end{abstract}
\section*{Introduction}
The study of higher order differential structures on Hopf 
algebras was initiated by S. L. Woronowicz \cite{Wo}, and 
further developed by K. Bresser {\it et al.}\cite{BMDS} and 
S. Majid \cite{Maj} for algebras of functions on finite 
groups. In particular, S. Majid has introduced and studied 
flat connections on the symmetric group $S_N.$ In our paper, 
we study the algebra generated by flat connections in a sense 
of Majid on a finite Weyl group. This is an interesting problem 
which is not treated in \cite{Maj}. 

We consider the differential structure with respect to 
the set of reflections. 
Since the complete set of the defining relations of the 
left-invariant exterior differential algebra has 
not yet been determined in general, 
we will work on its quadratic version $\Lambda_{quad}$ for the 
root system of type $A$ or $D,$ and on its quartic version 
$\Lambda_{quar}$ for the root system of type $B.$ 
Our main result describes a set of relations among flat 
connections on Weyl groups of type $A$ and $D.$ 
Conjecturally, 
these relations are complete set of relations among flat 
connections in $\Lambda_{quad}.$ 
We expect some connections of our construction with 
Schubert calculus on flag varieties \cite{KM}. 
\section{Woronowicz exterior algebra}
Woronowicz exterior algebra was introduced in \cite{Wo} 
for the study of higher order differential structure on the 
quantum groups. In the category of modules over a commutative 
algebra, the exterior products of a module are constructed 
by using the canonical action of the symmetric groups 
on the tensor products. In general, such a construction 
does not work in the category of bimodules over a noncommutative 
algebra because of lack of canonical action of the symmetric 
groups on the tensor products.  However, in the category 
of bimodules over a Hopf algebra, one can obtain a natural 
generalization of the exterior product. In this paper, all 
(Hopf) algebras are over a field $K$ of characteristic zero. 
Let $H$ be a Hopf algebra. 
\begin{dfn}
A bimodule $M$ over $H$ is called a left (resp. right) covariant 
bimolude if $M$ has a left (resp. right) $H$-comodule structure 
compatible with the bimodule structure. A bimodule $M$ is called 
a bicovariant bimodule if $M$ has left and right covariant 
bimodule structures and the left coaction and the right coaction 
commute. 
\end{dfn}
\begin{dfn} 
Let $M$ be a bicovariant bimodule over a Hopf algebra $H.$ 
An element $x\in M$ is called left (resp. right) invariant 
if $x$ is mapped to 
$1 \otimes x$ (resp. $x\otimes 1$) by the comodule 
action of $H.$  
\end{dfn}
\begin{lem}
There exists a unique braiding $\Psi : M \otimes_H M 
\rightarrow M \otimes _H M$ such that 
$\Psi(\omega \otimes \eta)=\eta \otimes \omega$ for 
left invariant $\omega$ and right invariant $\eta.$ 
\end{lem}
The homomorphism $\Psi$ induces a homomorphism 
$\Psi_i : M^{\otimes_H n} \rightarrow M^{\otimes_H n},$ 
$1\leq i \leq n,$ which acts as $\Psi$ on $i$-th and 
$(i+1)$-st components and acts identically on the other 
components. Take an element $w\in S_n$ and its reduced 
decomposition $w=s_{i_1}\cdots s_{i_l},$ $s_i= (i,i+1).$ 
Then we can associate a homomorphism 
$\Psi(w) : M^{\otimes_H n} \rightarrow M^{\otimes_H n}$ 
to the element $w$ by defining $\Psi(w)= \Psi_{i_1}\cdots 
\Psi_{i_l}.$ Since $\Psi_i$'s satisfy the braid relations, 
the homomorphism $\Psi(w)$ is independent of the choice of 
reduced decomposition of $w.$ 
Now we define the antisymmetrizer $A_n$ on $M^{\otimes_H n}$ 
by the formula 
\[ A_n = \sum_{w\in S_n} \textrm{sgn}(w) \Psi(w). \]  
\begin{dfn}
Woronowicz exterior algebra $\bigwedge M$ 
is a quotient of the tensor 
algebra of $M$ over $H$ by the kernel of the antisymmetrizer, 
i.e. {\rm 
\[ \bigwedge M := T_H M / \bigoplus_n \textrm{Ker} (A_n). \]} 
\end{dfn}
\section{Differential structure on the Weyl group}
\subsection{Differential structure on a finite group}
First of all, let us remind some fundamental 
facts on noncommutative differential structures on the 
finite group following \cite{BMDS} and \cite{Maj}. 
\begin{dfn} Let $A$ be a $K$-algebra. The 
first order differential structure of $A$ is a pair of 
$A$-bimodule $\Omega_A^1$ and a $K$-linear map $d: A 
\rightarrow \Omega_A^1$ such that the map $d$ satisfies 
the Leibniz rule $d(ab)=(da)b+a(db),$ for $a,b\in A,$ 
and the image of $d$ generates $\Omega_A^1$ as a left 
$A$-module. 
\end{dfn}
\begin{dfn}
Let $H$ be a Hopf algebra. The first order differential 
structure $(d: H \rightarrow \Omega_H^1)$ is said to be 
bicovariant if $\Omega_H^1$ has a structure of a bicovariant 
bimodule and the map $d$ is a bicomodule homomorphism. 
\end{dfn}
As a consequence of the construction in Section 1, we 
have the Woronowicz exterior algebra of 
a bicovariant differential structure  $\Omega_H^1$
of a Hopf algebra $H.$ 
\begin{dfn}
The Woronowicz exterior differential algebra $\Omega_w$ for 
the bicovariant differential structure of a Hopf algebra 
$H$ is a Woronowicz exterior algebra of $\Omega_H^1,$ i.e. 
$\Omega_w := \bigwedge \Omega_H^1.$ The left invariant 
subalgebra of $\Omega_w$ is denoted by $\Lambda_w.$ 
\end{dfn}
Let $G$ be a finite group and $H$ an algebra of 
functions on $G$ taking values on $K.$ 
Now we consider the 
differential structure on the Hopf algebra 
$H=K(G).$ The set of 
the delta functions $\{ \delta_g \; | \; g\in G \}$ can be 
taken as a linear basis of $H.$ 

Now we construct a canonical differential structure of the 
algebra $H.$ Take a subset ${\cal C}$ of $G$ 
which does not contain the identity element. Let 
$D_{\cal C}= \{ (x,y)\in G \times G \; | \; x^{-1}y 
\in {\cal C} \} .$ Define $\Omega^1(G)$ as a $K$-linear space 
generated by the set $\{ \delta_x \otimes \delta_y |(x,y)\in 
D_{\cal C} \} ,$ and 
\[ df= \sum_{(x,y)\in D_{\cal C}}(f(y)-f(x))\delta_x 
\otimes \delta_y , \; \; \; \; \textrm{for} \; \; 
f\in H. \]
Then $(d: H \rightarrow \Omega^1(G))$ is a first order 
differential structure on $H.$ All left covariant 
differential structures on $H$ are of this form, and 
$\Omega^1(G)$ is bicovariant if and only if the set 
${\cal C}$ is stable under the adjoint action of $G.$ 
Hence, simple bicovariant differential structures on 
$H$ are classified by nontrivial conjugacy classes of $G.$ 

For an element $a\in G,$ let $e_a=\sum_{g\in G} \delta_g 
d\delta_{ga}.$ Then the left invariant subalgebra 
$\Lambda_w$ is a $K$-linear subspace spanned by $e_a,$ 
$a\in {\cal C}.$ 
\subsection{Differential structure on the Weyl group}
Now we assume the group $W$ to be a Weyl group. Denote by 
$\Delta$ the set roots, and $\Delta_+$ the set of positive 
roots. As we have seen in the previous subsection, 
bicovariant differential structures on $H=K(W)$ 
are corresponding to adjoint invariant subsets of $W.$ 
We take ${\cal C}={\cal C}_{refl}$ the set of reflections 
as the simplest adjoint invariant subset of $W.$ 
\begin{rem}
{\rm For a simply-laced root system, the set ${\cal C}$ forms a 
conjugacy class. However, for a nonsimply-laced system 
the set ${\cal C}$ splits 
into a disjoint union of two conjugacy classes: 
${\cal C} = {\cal C}_l \cup {\cal C}_s,$ where 
${\cal C}_l$ (resp. ${\cal C}_s$) is the set of 
reflections with respect to the long (resp. short) roots. 
We can see that 
\[ \Lambda_w(B_n;{\cal C}_l) \cong  \Lambda_w(C_n;{\cal C}_s)
\cong \Lambda_w(D_n;{\cal C}), \] 
\[ \Lambda_w(B_n;{\cal C}_s) \cong  \Lambda_w(C_n;{\cal C}_l)
\cong \Lambda_w((A_1)^n;{\cal C}) , \] 
\[ \Lambda_w(G_2;{\cal C}_l) \cong  \Lambda_w(G_2;{\cal C}_s)
\cong  \Lambda_w(A_2;{\cal C}). \] 
This fact shows that the {\it simple} differential 
structure corresponding to ${\cal C}_l$ or ${\cal C}_s$ is 
not appropriate to investigate the differential structure 
for nonsimply-laced root systems. For that reason, we 
consider the differential structure obtained from the set 
${\cal C},$ which is not simple differential structure for 
nonsimply-laced root system. The algebra $\Lambda_w(X,{\cal C})$ 
will be denoted simply by $\Lambda_w(X).$ }
\end{rem}
We define a quadratic version of the left-invariant differential 
algebra follwing \cite{Maj}: 
\[ \Lambda_{quad} := T_K \Lambda^1 / \textrm{ker}(1- \Psi). \] 
\begin{conj}
{\rm  For simply-laced root systems, 
$\Lambda \cong \Lambda_{quad}.$ 
For the root system of type $A$ this conjecture was stated by 
S. Majid \cite{Maj}.} 
\end{conj}
\begin{rem} {\rm 
For nonsimply-laced root systems, $\Lambda_{quad}$ is not 
isomorphic to $\Lambda_w.$} 
\end{rem} 
\begin{ex}
{\rm The algebra $\Lambda_{quad}(B_n)$ is generated by 
$e_{(ij)},$ $e_{\overline{(ij)}}$ and $e_{(i)},$ where 
$(ij),$ $\overline{(ij)}$ and $(i)$ are reflections. 
The defining quadratic relations are: 
\[ e_{(ij)}^2=e_{\overline{(ij)}}^2=e_{(i)}^2=0, \] 
\[ e_{(ij)}e_{(kl)}+e_{(kl)}e_{(ij)}= 
e_{(ij)}e_{\overline{(kl)}}+e_{\overline{(kl)}}e_{(ij)}
=e_{\overline{(ij)}}e_{\overline{(kl)}}+
e_{\overline{(kl)}}e_{\overline{(ij)}} =0, \; \; \; 
\textrm{for} \; \; \{ i,j \} \cap \{ k,l \} = \emptyset , \] 
\[ e_{(i)}e_{(j)}+e_{(j)}e_{(i)}= 
e_{(ij)}e_{\overline{(ij)}}+e_{\overline{(ij)}}e_{(ij)}=
e_{(ij)}e_{(k)}+e_{(k)}e_{(ij)}=
e_{\overline{(ij)}}e_{(k)}+e_{(k)}e_{\overline{(ij)}}, \; \; 
\textrm{if} \; \; k\not= i,j, \]
\[ e_{(ij)}e_{(jk)}+e_{(jk)}e_{(ki)} + e_{(ki)}e_{(ij)}=0, \] 
\[  e_{\overline{(ik)}}e_{(ij)}+e_{(ji)}e_{\overline{(jk)}} + 
e_{\overline{(kj)}}e_{\overline{(ik)}}=0, \] 
\[ e_{(ij)}e_{(i)}+e_{(j)}e_{(ij)}+
e_{(i)}e_{\overline{(ij)}}+e_{\overline{(ij)}}e_{(j)}=0. \] 
The algebra $\Lambda_{quad}(D_n)$ is a quotient of 
$\Lambda_{quad}(B_n).$}
\end{ex}
\begin{rem} {\rm 
The algebra $\Lambda_{quad}(B_n)$ is not isomorphic to 
$\Lambda_w(B_n).$ For example, the relations 
\[ e_{\overline{(ij)}}e_{(i)}e_{(ij)}e_{(i)}+ 
e_{(i)}e_{(ij)}e_{(i)}e_{\overline{(ij)}}+
e_{(ij)}e_{(i)}e_{\overline{(ij)}}e_{(i)}+
e_{(i)}e_{\overline{(ij)}}e_{(i)}e_{(ij)} =0, \] 
\[ e_{(ij)}e_{(i)}e_{(ij)}e_{(i)}+ 
e_{(i)}e_{(ij)}e_{(i)}e_{(ij)} =0 \] 
hold in $\Lambda_w(B_n),$ but they do not in 
$\Lambda_{quad}(B_n).$ We denote by $\Lambda_{quar}(B_n)$ 
the quotient algebra of $\Lambda_{quad}(B_n)$ by the ideal 
generated by the quartic relations above.} 
\end{rem}
\begin{ex} {\rm 
The algebra $\Lambda_{quad}(B_2)$ is infinite dimensional. It 
has the Hilbert polynomial 
\[ 1+4t+8t^2+12t^3+16t^4+20t^5+24t^6+28t^7+32t^8+\cdots =
(1+t)^2(1-t)^{-2}. \]
In the algebra $\Lambda_w(B_2),$ the quartic relations 
\[ e_{\overline{(12)}}e_{(1)}e_{(12)}e_{(1)}+ 
e_{(1)}e_{(12)}e_{(1)}e_{\overline{(12)}}+
e_{(12)}e_{(1)}e_{\overline{(12)}}e_{(1)}+
e_{(1)}e_{\overline{(12)}}e_{(1)}e_{(12)} =0 \] 
and 
\[ e_{(12)}e_{(1)}e_{(12)}e_{(1)}+ 
e_{(1)}e_{(12)}e_{(1)}e_{(12)} =0 \] 
hold. The algebra $\Lambda_{quar}(B_2)$ obtained by adding the 
quartic relations above 
to $\Lambda_{quad}(B_2)$ is finite dimensional and has the 
Hilbert polynomial 
\[ (1+t)^4(1+t^2)^2. \] 
In particular, $\Lambda_w(B_2)$ is finite dimensional. 
The anticommutative quotient of the algebra 
$\Lambda_{quad}(B_2)$ has the Hilbert polynomial 
\[ 1+4t+5t^2+2t^3=(1+t)^2(1+2t). \] }
\end{ex}
\begin{conj} {\rm 
The relations in Example 2.1 and Remark 2.3 are the complete 
set of relations for $\Lambda_w(B_n),$ i.e. 
$\Lambda_w(B_n) \cong \Lambda_{quar}(B_n).$ } 
\end{conj}
\section{$U(1)$-gauge theory}
The algebra $\Lambda_w$ has a structure of a differential graded 
algebra over $K.$ Denote by $H^*(W)$ the cohomology group of the differential 
graded algebra $\Lambda_w.$ 
Let $\theta=\sum_{a\in {\cal C}}e_a.$ If $W$ is nonsimply-laced, 
we define $\theta_1=\sum_{a\in {\cal C}_l}e_a$ and 
$\theta_2=\sum_{a\in {\cal C}_s}e_a.$ 
\begin{prop}
For simply-laced root system, 
$H^1(W)= K\cdot \theta.$ For nonsimply-laced root system, 
$H^1(W)= K\cdot \theta_1 \oplus K\cdot \theta_2.$ 
\end{prop}
{\it Proof} \quad Since $d\delta_g = \sum_{c\in {\cal C}} 
(\delta_{gc} - \delta_g)e_c$ and $e_c \cdot \delta_g = \delta_{gc} 
\cdot e_c,$ we have 
\[ de_a = \sum_{g\in W} d\delta_g d\delta_{ga} = \theta e_a + e_a 
\theta. \] 
We can show that if $\eta =\sum_{a\in {\cal C}} \eta_a e_a$ 
is a closed 1-form, then $\eta_a = \eta_b$ must be satisfied 
when $a$ and $b$ are conjugate each other. \quad \rule{3mm}{3mm} 

Let $\Omega^1_H$ be a bicovariant differential structure 
of a Hopf algebra $H.$ 
\begin{dfn}
For a 1-form $\eta \in \Omega^1_H,$ the covariant 
curvature is defined by 
\[ F(\eta)= d \eta + \eta \wedge \eta . \] 
If $F(\eta)= 0,$ $\eta$ is called a flat $U(1)$-connection. 
\end{dfn} 
As we have seen in Remark 2.1, the simple differential 
structures of nonsimply-laced Weyl groups are reduced to 
the ones of simply-laced Weyl groups. Hence, we restrict 
our considerations to the case of simply-laced root system 
$X=(\Delta, W).$ Fix the set of simple roots $\Sigma \subset 
\Delta.$ 
Let $\omega_{\alpha}$ be a fundamental dominant dominant 
weight corresponding to a simple root $\alpha \in \Sigma.$
Denote by $(\nu_{\alpha})_{\alpha}$ Schmidt's orthogonalization of 
$(\omega_{\alpha})_{\alpha}$. 
We define the 1-forms $\theta_{\alpha}^X$ for $\alpha \in \Sigma$ 
by 
\[ \theta_{\alpha}^X:= \sum_{\gamma \in \Delta_+(\alpha)} \langle
\nu_{\alpha}, \gamma^{\vee} \rangle e_{s_{\gamma}} ,\] 
where $\Delta_+(\alpha)$ is the set of the roots $\gamma$ satisfying 
the condition $\langle \nu_{\alpha} , \gamma \rangle >0.$ 
\begin{prop}
For the classical root systems the 1-forms 
$\theta_{\alpha}=\theta_{\alpha}^X$ satisfy 
the relations of anticommutativity 
$\theta_{\alpha}\theta_{\beta}+\theta_{\beta}\theta_{\alpha}=0$ 
and flatness relations $F(-\theta_{\alpha})=0.$  
\end{prop}
{\it Proof} \quad This can be shown by direct computations. (See Section 
5.) \quad \rule{3mm}{3mm} 
\begin{rem}{\rm 
The 1-forms $\theta_{\alpha}$ can be considered as an 
analogue of the Dunkl elements introduced in \cite{FK} and \cite{KM}. 
}
\end{rem}
\begin{ex}
{\rm Here we give an example in the exceptional root system of 
type $G_2.$ Let $\alpha$ be the short simple root and $\beta$ the long 
one. Then the set of positive roots is 
\[ \Delta_+= \{ a_1= \alpha, a_2= 3\alpha+ \beta, a_3= 3\alpha+ 2\beta, 
a_4= 2\alpha+\beta, a_5=\alpha+\beta,a_6=\beta \} . \] 
Let $s_i$ be the reflection with respect to $a_i$ and 
$e_i:=e_{s_i}.$ Then the relations 
\[ e_i^2=0, \; \; e_1e_4+e_4e_1=e_2e_5+e_5e_2=e_3e_6+e_6e_3=0, \]
\[ e_1e_3+e_3e_5+e_5e_1=e_3e_1+e_5e_3+e_1e_5=0, \] 
\[ e_2e_4+e_4e_6+e_6e_2=e_4e_2+e_6e_4+e_2e_6=0, \] 
\[ e_1e_2+e_2e_3+e_3e_4+e_4e_5+e_5e_6+e_6e_1=0, \] 
\[ e_2e_1+e_3e_2+e_4e_3+e_5e_4+e_6e_5+e_1e_6=0 \] 
hold in $\Lambda_w(G_2).$ The first cohomology group is 
\[ H^1(W(G_2))= K\cdot (e_1+e_3+e_5)+K\cdot (e_2+e_4+e_6). \] 
Moreover, the 1-forms $\eta_1=-(2e_1+e_2+e_3+e_5+e_6)$ and 
$\eta_2=-(e_2+e_3+2e_4+e_5+e_6)$ define flat connections which 
satisfy the anticommutativity $\eta_1 \eta_2+\eta_2 \eta_1=0.$ 
}
\end{ex}
\section{Hopf algebra structure}
We introduce a Hopf algebra structure on 
$K\langle W \rangle \otimes_K \Lambda_{quad}(W).$ 
We consider $K\langle W \rangle \otimes_K \Lambda_{quad}(W)$ 
as a twisted group algebra defined by the commutation relations 
\[ e_{s_{\gamma}}\cdot w = (-1)^{l(w)}w\cdot e_{s_{w\gamma}}, \]
where $s_{\gamma}$ is a reflection with respect to a root 
$\gamma$ and $w\in W.$ The coproduct $\Delta,$ the antipode 
$S$ and the counit $\varepsilon$ are given by the formulas: 
\[ 
\begin{array}{ll}
\Delta(e_{s_{\gamma }}) = e_{s_{\gamma}} \otimes 1 + s_{\gamma} 
\otimes e_{s_{\gamma}}, & \Delta (w) = w \otimes w , \\
S( e_{s_{\gamma}} ) = -s_{\gamma} \cdot e_{s_{\gamma}}, & 
S(w)= w^{-1}, \\  
\varepsilon(e_{s_{\gamma}}) =0, & \varepsilon (w)=1. 
\end{array}
\] 
The adjoint representation of the Hopf algebra gives an 
action of $\Lambda_{quad}(W)$ on itself. The element 
$e_{s_{\gamma}}$ acts as a twisted derivation 
\[ D_{\gamma} (x) = e_{s_{\gamma}}x-(-1)^{\deg x}s_{\gamma}(x)
e_{s_{\gamma}}, \] 
for a homogeneous element $x\in \Lambda_{quad}(W).$ 
The twisted derivation $D_{\gamma}$ satisfies 
the twisted Leibniz rule 
\[ D_{\gamma} (xy)= D_{\gamma}(x)y+(-1)^{\deg x}s_{\gamma}(x)
D_{\gamma}(y). \] 
\begin{rem}
{\rm The Hopf algebra considered above coincides with 
the one obtained as a twisted group algebra over the quadratic 
lift of the bracket algebra $BE(W,S)$ defined in \cite{KM}. 
(If the root system is simply-laced, the bracket algebra 
itself is a quadratic algebra.) 
In particular, it coincides with the fibered Hopf algebra 
introduced in \cite{FP} for the root system of type $A.$ }
\end{rem}
\section{Subalgebra generated by flat connections}
In this section, we discuss on the structure of the subalgebra 
generated by the flat connections $\theta_{\alpha},$ which are 
introduced in Section 3. We will treat only classical 
root systems. Since we use only quadratic relations, 
we work on the quadratic algebra $\Lambda_{quad}.$ 
For simplicity, the symbols $(ij),$ $\overline{(ij)}$ and 
$(i)$ are used instead of $e_{(ij)},$ $e_{\overline{(ij)}}$ 
and $e_{(i)}$ respectively. 
The $U(1)$-connections 
$\theta_1^X,\ldots,\theta_n^X \in \Lambda_{quad}(X)$ 
$(X=A_{n-1},B_n,D_n)$ 
are defined as follows: 
\[ \theta^{A_{n-1}}_i= \sum_{j=1}^n (ij), \]
\[ \theta^{D_n}_i=\sum_{j=1}^n ((ij)+\overline{(ij)}), \] 
\[ \theta^{B_n}_i=\sum_{j=1}^n ((ij)+\overline{(ij)})+
2(i). \]
We can easily check that the elements $-\theta_i$ define flat connections 
and satisfy the anticommutativity 
$\theta_i\theta_j+\theta_j\theta_i=0.$ For example, 
\[ \theta_i^{B_n}\theta_j^{B_n}+\theta_j^{B_n}\theta_i^{B_n} = \]
\[ \sum_{k\not= i,j} \sum_{\varsigma_1,\varsigma_2=\pm}
\Bigl((ij)^{\varsigma_1}(jk)^{\varsigma_2}+
(ik)^{\varsigma_1}(ij)^{\varsigma_2}+(ik)^{\varsigma_1}(jk)^{\varsigma_2}+
(jk)^{\varsigma_1}(ij)^{\varsigma_2}+(ij)^{\varsigma_1}(ik)^{\varsigma_2}+
(jk)^{\varsigma_1}(ik)^{\varsigma_2} \Bigr) \] 
\[ + 2\sum_{k\not= i,j}\sum_{\varsigma=\pm}
\Bigl( (ik)^{\varsigma}(j)+(j)(ik)^{\varsigma}+
(jk)^{\varsigma}(i)+(i)(jk)^{\varsigma} \Bigr) =0, \]
where $(ij)^+=(ij)$ and $(ij)^-=\overline{(ij)}.$ Since $\sum_i \theta_i^{B_n} 
=2\theta^{B_n},$ the flatness $-d\theta_i^{B_n}+(\theta_i^{B_n})^2=0$ follows 
from the anticommutativity. 
\begin{lem}
{\rm (Cyclic relations in $\Lambda_{quad}(A_{n-1})$)} \\ 
For any distinct $1\leq a_1,\ldots,a_k \leq n,$ 
\[ \sum_{i=2}^k (-1)^{k(i-1)} (a_1a_i)(a_1a_{i+1})\cdots 
(a_1a_k)(a_1a_2)\cdots (a_1a_i) =0. \]
\end{lem}
{\it Proof} \quad These relations are obtained by 
applying the composition of twisted derivations 
$D_{a_{k-1}a_k}D_{a_{k-2}a_{k-1}}\cdots D_{a_2a_3}$ to the 
relation $(a_1a_2)^2=0.$ \quad \rule{3mm}{3mm}
\begin{ex} {\rm 
Let $a_1,a_2,a_3,a_4$ be distinct, then 
\begin{itemize}
\item $k=3$ \hspace{10mm} 
$(a_1a_2)(a_1a_3)(a_1a_2)-(a_1a_3)(a_1a_2)(a_1a_3)=0$ 
\item $k=4$ 
\[ (a_1a_2)(a_1a_3)(a_1a_4)(a_1a_2)+(a_1a_3)(a_1a_4)(a_1a_2)
(a_1a_3)+(a_1a_4)(a_1a_2)(a_1a_3)(a_1a_4)=0 \]
\end{itemize}}
\end{ex}
\begin{lem}
For any distinct $1\leq a_1,\ldots,a_{k+1} \leq n,$ 
\[ \Bigl( \prod_{j=2}^k (a_1a_j) \Bigr) (a_1a_2)(a_1a_{k+1}) + 
(-1)^{k+1} (a_1a_{k+1}) 
\Bigl( \prod_{j=2}^k (a_1a_j) \Bigr) (a_1a_2) \] 
\[ +\Bigl( \prod_{j=2}^{k+1} (a_1a_j) \Bigr) (a_2a_{k+1})+ 
(-1)^{k+1}(a_2a_{k+1})(a_1a_{k+1})
\Bigl( \prod_{j=3}^k (a_1a_j) \Bigr) (a_1a_2) =0. \]
\end{lem}
{\it Proof} \quad By using the equalities 
$(a_1a_{k+1})(a_2a_{k+1})+(a_2a_{k+1})(a_1a_2)+
(a_1a_2)(a_1a_{k+1})=0,$ 
$(a_2a_{k+1})(a_1a_{k+1})+(a_1a_{k+1})(a_1a_2)+
(a_1a_2)(a_2a_{k+1})=0$ 
and anticommutativity relations, we obtain 
\[ \Bigl( \prod_{j=2}^k (a_1a_j) \Bigr) (a_1a_2)(a_1a_{k+1})+
\Bigl( \prod_{j=2}^{k+1} (a_1a_j) \Bigr) (a_2a_{k+1}) = -
\Bigl( \prod_{j=2}^k (a_1a_j) \Bigr) (a_2a_{k+1})(a_1a_2)  \]
\[ =  -(-1)^{k-2}(a_1a_2)(a_2a_{k+1})
\Bigl( \prod_{j=3}^k (a_1a_j) \Bigr) (a_1a_2)  \] 
\[ =  (-1)^k (a_1a_{k+1}) 
\Bigl( \prod_{j=2}^k (a_1a_j) \Bigr) (a_1a_2) + 
(-1)^k (a_2a_{k+1})(a_1a_{k+1})
\Bigl( \prod_{j=3}^k (a_1a_j) \Bigr) (a_1a_2).  \]
This completes the proof. \quad \rule{3mm}{3mm} 
\begin{cor}
For any $m\in {\bf Z}_{\geq 1},$ 
\[ (\theta_1^{A_{n-1}})^{2m}+\cdots +(\theta_n^{A_{n-1}})^{2m}=
0. \]
\end{cor}
{\it Proof}\quad The sum $\sum_i \theta_i^{2m}$ is a 
sum of products of cycles, and the number of odd cycles is even. 
All even cycles give zero contribution, see Lemma 5.1. 
According to Lemma 5.2 we can kill all even products of 
odd cycles. \quad \rule{3mm}{3mm}
\begin{lem} 
\abovedisplayskip -2mm 
\[ \theta_1^{A_{n-1}}\cdots \theta_n^{A_{n-1}} =0. \] 
\end{lem}
\begin{lem} For any integer $k$ between $1$ and $n$, we have
$$  \prod_{j=1, j \ne k}^{n} \theta_j^{A_{n-1}}= \sum_{\sigma 
\in Per(1,\cdots,\widehat k,\cdots,n)} 
(-1)^{l(\sigma)}\prod_{j=1, j \ne k}^{n}(\sigma(j),k),$$
where the sum runs over all permutations $\sigma$ of the set $(1,\cdots,
\widehat k,\cdots n)$ and $l(\sigma)$ denotes the length of permutation $\sigma.$
\end{lem}
{\it Proof of Lemmas 5.3 and 5.4} \quad The proof is by induction 
on $n.$ We will prove the equations in Lemmas 5.3 and 5.4 for $A_{n-1}$ 
under the assumption that Lemma 5.3 holds for $A_{n-2}.$ The equation 
$\theta_1^{A_{n-2}}\cdots \theta_{n-1}^{A_{n-2}} =0$ means that we have in 
$\Lambda_{quad}(A_{n-1})$ 
\[ \sum_{i_1,\ldots ,i_{n-1}} (1,i_1)\cdots (n-1,i_{n-1}) =0, \] 
where $i_1,\ldots ,i_{n-1}$ run over the letters satisfying $i_l\not= l$ and 
$1\leq i_1,\ldots ,i_{n-1} \leq n-1.$ Let 
$M_1$ be the sum of the products $\prod_{j=1, j \ne k}^n (j,i_j)$ such 
that none of the letters $i_j$ equal $k,$ and $M_2$ be the sum of the products 
such that at least one letter $i_j$ equals $k.$ Then, 
$\prod_{j=1, j \ne k}^{n} \theta_j^{A_{n-1}}=M_1+M_2.$ 
The assumption of the induction shows $M_1=0.$ 
We can express $\prod_{j=1, j \ne k}^{n}(\sigma(j),k)$ as a sum of 
terms of form $\pm (1,b_1)\cdots (n, b_n)$ by applying substitution 
$(a_i b)(a_{i+1} b)\rightarrow -(a_{i+1} b)(a_i a_{i+1})-
(a_{i+1}a_i)(a_ib)$ repeatedly 
when a term $\cdots (a_i b)(a_{i+1} b) \cdots$ 
with $a_i>a_{i+1}$ appears. This procedure yields the equality 
\[ \sum_{\sigma 
\in Per(1,\cdots,\widehat k,\cdots,n)} 
(-1)^{l(\sigma)}\prod_{j=1, j \ne k}^{n}(\sigma(j),k) = M_2. \] 
Now we have the equality in Lemma 5.4. Multiply both hand side 
by $\theta_k^{A_{n-1}}.$ Then we have 
\[ (-1)^{k-1} \theta_1^{A_{n-1}}\cdots \theta_n^{A_{n-1}} = 
\sum_{\sigma 
\in Per(1,\cdots,\widehat k,\cdots,n)} \sum_{l\ne k}
(-1)^{l(\sigma)}(k,l) \prod_{j=1, j \ne k}^{n}(\sigma(j),k). \] 
Here, we can show that the right hand side is equal to zero from 
the cyclic relations in Lemma 5.1. \quad \rule{3mm}{3mm} 
\begin{lem}
$$\sum_{k=1}^{m}(-1)^{(m-1)(k-1)}\prod_{j=k+1}^{m}(k,j)
\prod_{j=1}^{k-1}
(j,k)=0.$$
\end{lem}
{\it Proof} \quad By induction, one can show 
\[ D_{a_ma_{m+1}}\Bigl( \prod_{j=1}^{m-1}(a_ja_m) \Bigr) 
= (-1)^m(a_ma_{m+1})\Bigl( \prod_{j=1}^{m-1} (a_j a_m) \Bigr)
+ \prod_{j=1}^m (a_j a_{m+1}) \] 
by using the identity 
\[  \prod_{j=1}^{m-1}(a_ja_m)= (-1)^{n-2} \prod_{j=2}^{m-1}
(a_ja_m)\cdot (a_1a_2) - (a_1a_2)(a_1a_m)\prod_{j=3}^{m-1}
(a_ja_m). \] 
Then, the desired identity is obtained by applying 
$D_{a_{m-1}a_m}\cdots D_{a_3a_4}$ to the identity 
$(a_1a_2)(a_1a_3)+(a_2a_3)(a_1a_2)+(a_1a_3)(a_2a_3)=0.$
\quad \rule{3mm}{3mm}
\begin{ex} {\rm 
Let $a_1,a_2,a_3,a_4$ be distinct, then 
\begin{itemize}
\item $m=3$ \hspace{10mm} 
$(a_1a_2)(a_1a_3)+(a_2a_3)(a_1a_2)+(a_1a_3)(a_2a_3)=0.$ 
\item $m=4$ 
\[ (a_1a_2)(a_1a_3)(a_1a_4)-(a_2a_3)(a_2a_4)(a_1a_2)
+(a_3a_4)(a_1a_3)(a_2a_3)-(a_1a_4)(a_2a_4)(a_3a_4)=0. \]
\end{itemize}}
\end{ex}
\begin{thm}
The connections $\theta_1^{A_{n-1}},\ldots,
\theta_n^{A_{n-1}}$ satisfy 
the following relations: 
\[ \epsilon_k((\theta_1^{A_{n-1}})^2,\ldots,(\theta_n^{A_{n-1}})^2)=0, 
\; \; \; \; 1\leq k \leq n, \] 
where $\epsilon_k$ is the $k$-th elementary symmetric polynomial. 
Moreover, 
\[ \theta_1^{A_{n-1}}\cdots \theta_n^{A_{n-1}}=0, \]
\[ \sum_{i=1}^n (-1)^{i+1}\theta_1^{A_{n-1}}\cdots 
\hat{\theta}_i^{A_{n-1}}\cdots \theta_n^{A_{n-1}} =0. \]
\end{thm}
{\it Proof}\quad Indeed, the first series of equalities 
follow from Corollary 5.1. The second equality has been proved in 
Lemma 5.3. 
The last relation follows from 
Lemmas 5.4 and 5.5. 
\quad \rule{3mm}{3mm} 
\medskip \\ 
Let us remark that 
\[ \epsilon_{n-1}((\theta_1^{A_{n-1}})^2,\ldots,(\theta_n^{A_{n-1}})^2)
= \left( \sum_{i=1}^n (-1)^{i+1}\theta_1^{A_{n-1}}\cdots 
\hat{\theta}_i^{A_{n-1}}\cdots \theta_n^{A_{n-1}} \right) ^2. \] 
\begin{prop}
The elements $E_{(ij)}:=e_{(ij)}+e_{\overline{(ij)}}\in 
\Lambda_{quad}(D_n)$ generate a subalgebra isomorphic to 
$\Lambda_{quad}(A_{n-1}),$ where we have the natural 
identification $\theta_j^{A_{n-1}}
=\theta_j^{D_n},$ $1\leq j \leq n.$ 
\end{prop}
{\it Proof} \quad We can check the identities 
\[ E_{(ij)}^2=0, \] 
\[ E_{(ij)}E_{(kl)}+E_{(kl)}E_{(ij)}=0, \; \; \; 
\textrm{for} \; \; \{ i,j \} \cap \{ k,l \} = \emptyset , \] 
\[ E_{(ij)}E_{(jk)}+E_{(jk)}E_{(ki)} + E_{(ki)}E_{(ij)}=0. \] 
Hence, we can define an algebra homomorphism 
$\iota : \Lambda_{quad}(A_{n-1}) \rightarrow \Lambda_{quad}(D_n)$ 
by mapping $e_{(ij)}$ to $E_{(ij)}.$ We also have an algebra 
homomorphism $\pi : \Lambda_{quad}(D_n) \rightarrow 
\Lambda_{quad}(A_{n-1})$ 
obtained by putting $e_{\overline{(ij)}}=0.$ Since 
$\pi \circ \iota = {\rm id},$ 
the elements $E_{(ij)}$ generate a subalgebra isomorphic to 
$\Lambda_{quad}(A_{n-1}).$ \quad \rule{3mm}{3mm}
\begin{cor}
\[ \epsilon_k((\theta_1^{D_n})^2,\ldots,(\theta_n^{D_n})^2)=0, 
\; \; \; \; 1\leq k \leq n. \] 
Moreover, 
\[ \theta_1^{D_n}\cdots \theta_n^{D_n}=0, \]
\[ \sum_{i=1}^n (-1)^{i+1}\theta_1^{D_n}\cdots 
\hat{\theta}_i^{D_n}\cdots \theta_n^{D_n} =0. \]
\end{cor}
\begin{conj}{\rm (1) Let $X$ denote either $A_{n-1}$ or 
$D_n.$ Relations 
\[ \epsilon_k((\theta_1^X)^2,\ldots,(\theta_n^X)^2)=0, 
\; \; \; \; 1\leq k \leq n. \] 
\[ \theta_1^X\cdots \theta_n^X=0, \]
\[ \sum_{i=1}^n (-1)^{i+1}\theta_1^X\cdots 
\hat{\theta}_i^X\cdots \theta_n^X =0 \]
together 
with the anticommutativity relations $\theta_i^X\theta_j^X+
\theta_j^X\theta_i^X=0,$ form the complete list of relations 
among $\theta_1^X,\ldots ,\theta_n^X$ in the quadratic 
algebra $\Lambda_{quad}(X).$ \medskip \\ 
(2) For $X=B_n,$ the relations 
\[ \epsilon_k((\theta_1^X)^2,\ldots,(\theta_n^X)^2)=0, 
\; \; \; \; 1\leq k \leq n \] 
and the anticommutativity relations form the complete list of relations 
among $\theta_1^{B_n},\ldots ,\theta_n^{B_n}$ in the algebra 
$\Lambda_{quar}(B_n).$}
\end{conj}
We can check that the above relations are 
valid in the algebra $\Lambda_{quar}(B_n)$ for $n\leq 3.$ 
\begin{rem}
{\rm Let us consider the flag variety $Fl_n$ of type $A_{n-1}$ 
and the tautological flag on it: 
\[ 0=F_0 \subset F_1 \subset F_2 \subset \cdots \subset F_n = 
{\cal O}_{Fl_n}^{\oplus n}. \] 
The cohomology ring $H^*(Fl_n,K)$ is isomorphic to the algebra 
\[ K[x_1,\ldots x_n]/(\epsilon_1(x),\ldots ,\epsilon_n(x)), \] 
where $x_i=c_1(F_i/F_{i-1}).$ The algebra generated by the flat 
connections $\theta_i^{A_{n-1}}$ can be considered as a super-analogue 
of the cohomology ring of the flag variety, and our result shows 
that both algebras have some common relations in even degrees. }
\end{rem}

Research Institute for Mathematical Sciences \\ 
Kyoto University \\
Sakyo-ku, Kyoto 606-8502, Japan \\
e-mail: kirillov@kurims.kyoto-u.ac.jp 
\bigskip \\ 
Department of Mathematics \\ 
Kyoto University \\ 
Sakyo-ku, Kyoto 606-8502, Japan \\ 
e-mail: maeno@math.kyoto-u.ac.jp 
\end{document}